\documentclass[11pt,reqno,a4paper]{amsart}%
\usepackage{amsfonts,amsmath,times, amssymb}
\usepackage{a4wide}
\usepackage{tikz}
\definecolor{blau}{rgb}{0,0,0.75} 
\usepackage[pdftex]{hyperref}
\hypersetup{colorlinks,linkcolor=blau,citecolor=blue,urlcolor=blau}

\allowdisplaybreaks

\newcommand{\fallfak}[2]{\ensuremath{#1^{\underline{#2}}}}
\newcommand{\Stir}[2]{\genfrac{ \{ }{ \} }{0pt}{}{#1}{#2}}
\newcommand{\stir}[2]{\genfrac{ [ }{ ] }{0pt}{}{#1}{#2}}
\newcommand{\N}{\ensuremath{\mathbb{N}}}
\newcommand{\law}{\ensuremath{\stackrel{(d)}=}}

\def\P{{\mathbb {P}}}
\def\E{{\mathbb {E}}}
\def\V{{\mathbb {V}}}

\newtheorem{theorem}{Theorem}
\newtheorem{lemma}{Lemma}
\newtheorem{prop}{Proposition}
\newtheorem{coroll}{Corollary}

\theoremstyle{definition}
\newtheorem{remark}{Remark}

\title{On generalized P\'olya urn models}

\author[M.-R.~Chen]{May-Ru Chen}
\address{May-Ru Chen\\
Department of Applied Mathematics\\
National Sun Yat-sen University\\
No. 70, Lienhai Road\\
Kaohsiung 80424, Taiwan, R.O.C. }
\email{chenmr@math.nsysu.edu.tw}

\author[M.~Kuba]{Markus Kuba}
\address{Markus Kuba\\
Institut f{\"u}r Dis\-krete Mathe\-matik und Geo\-metrie\\
Technische Uni\-versit\"at Wien\\
Wiedner Hauptstr. 8-10/104\\
1040 Wien, Austria -- HTL Wien 5 Spengergasse, Spengergasse 20, 1050 Wien, Austria} %
\email{kuba@dmg.tuwien.ac.at}

\thanks{The second author was partially supported by the Austrian Science Foundation FWF, grant S9608-N13.}

\date{\today}


\begin{document}

\begin{abstract}
We study an urn model introduced in the paper of Chen and Wei~\cite{ChenWei}, where at each discrete time step $m$ balls are drawn at random 
from the urn containing colors white and black. Balls are added to the urn according to the inspected colors, generalizing the well known P\'olya-Eggenberger urn model, case $m=1$.
We provide exact expressions for the expectation and the variance of the number of white balls after $n$ draws, 
and determine the structure of higher moments. Moreover, we discuss extensions to more than two colors. Furthermore, we introduce and discuss a new urn model 
where the sampling of the $m$ balls is carried out in a step-by-step fashion, and also introduce a generalized Friedman's urn model.
\end{abstract}

\maketitle

\section{Introduction}
P\'olya-Eggenberger urn models are defined as follows. At the start, time zero, 
the urn contains $W_0$ white balls and $B_0$ black balls. The
evolution of the urn occurs in discrete time steps. At every step
a ball is chosen at random from the urn. The color of the ball is
inspected and then the ball is reinserted into the urn. According
to the observed color of the ball, balls are added/removed due to
the following rules. If we have chosen a white ball, we put into
the urn $a$ white balls and $b$ black balls, but if we have chosen
a black ball, we put into the urn $c$ white balls and $d$ black
balls. The values $a, b, c, d \in \mathbb{Z}$ are fixed integer
values and the urn model is specified by the $2\times 2$ ball
replacement matrix $M = \bigl(\begin{smallmatrix} a & b \\ c &
d\end{smallmatrix}\bigr)$. Urn models are simple, useful mathematical tools for describing many
evolutionary processes in diverse fields of application such as
analysis of algorithms and data structures, statistics and genetics, see Johnson and Kotz~\cite{JohnsonKotz1977}, Kotz and N.~Balakrishnan~\cite{Kotz1997}, and also Mahmoud~\cite{Mah2009}.

\smallskip

One of the most fundamental urn models is original the P\'olya-Eggenberger urn model~\cite{Eggenberger1923}, associated with the ball
replacement matrix  $M = \left(\begin{smallmatrix} 1 & 0 \\ 0 &
1\end{smallmatrix}\right)$.  The P\'olya-Eggenberger urn is a balanced urn model, the total number of added or removed balls is
constant, independently of the observed color. A parameter of interest is
the number $W_n$ of white balls contained in the urn after $n$ draws. 
Various generalizations of this urn model have been considered by Hill~\cite{Hill1980}, Bagci and Pal~\cite{Bagci1985}, Pemantle~\cite{Pemantle1990}, Gouet~\cite{Gouet1989,Gouet1993}, Schreiber~\cite{Schreiber2001}; we refer to Chen and Wei~\cite{ChenWei} for a brief discussion of the aforehand mentioned works.

\smallskip

This work is devoted to the study of a generalization of the P\'olya-Eggenberger urn model, where \emph{several balls} are drawn at each discrete time step, their colors are inspected, and they reinserted into the urn. The addition/removal of balls depends on the combinations of colors of the drawn balls. Such urn models recently received some  attention in the literature, see for example Chen and Wei~\cite{ChenWei}, Mahmoud~\cite{Mah2009}, and Renlund~\cite{Renlund}. Chen and Wei~\cite{ChenWei} introduced a particular urn model they called model $M$, where $m\ge 1$ balls are drawn from the urn at each discrete time step. 
Say $m - \ell$ white balls and $\ell$ black balls have been drawn, $0\le\ell\le m$, their colors are noted, and the drawn balls are returned to the
urn together with addition $c(m - \ell)$ white balls and $c\ell$ black balls. The ball replacement matrix of this urn model is a rectangular matrix $M$, given by 

\begin{equation}
\label{NEWmodelM}
	M=\left[
	\begin{array}{cc}
	mc&	0	\\
	(m-1)c	& c\\
	\hdots&\hdots\\
	c	& (m-1)c\\
	0&	mc	\\
	\end{array}
	\right],
\end{equation}

with parameter $c\in\N$ and $m\ge 1$. The rows of the rectangular replacement matrix encode the sampling scheme in the obvious way, 
the $\ell$-th row corresponds to the case of drawing a combination of $(m-\ell)$ white balls and $\ell$ black balls, $0\le \ell \le m$, 
where $c(m - \ell)$ white balls and $c\ell$ black balls are being added. Note that in model $M$ the drawing of the $m$ balls occurs without replacement, in other words the distribution of the number of white balls in the sample of size $m$ follows a hypergeometric distribution. 

\smallskip

Chen and Wei studied the distribution of the number of white balls $W_n$ after $n$ draws and could show the almost sure convergence of $W_n$, suitably normalized,
to a continuous distribution by using martingales. The aim of this note is to provide further insight into the limiting distribution of 
the number of white balls by providing exact expressions for the expectation and the variance of $W_n$, from which one obtains the expectation and variance of the limit law. Moreover, we also show to obtain in principle exact expressions for arbitrary moments of the limit law. Note that the case $m=1$ corresponds to the well known P\'olya-Eggenberger urn $M = \left(\begin{smallmatrix} c & 0 \\ 0 &
c\end{smallmatrix}\right)$, which is completely understood, using different arguments such as counting arguments, or stochastic processes -- see for example Janson~\cite{Jan2006}.
It is known that proportion of white balls after $n$ draws is a martingale, and has a beta distribution as the limit law with parameters $b/c$ and
$w/c$. Therefore, the case $m=1$ is excluded from our study.

\smallskip 

Throughout this work we use the notations $\stir{n}k$ and $\Stir{n}k$, where $\stir{n}{k}$ denotes the unsigned Stirling numbers of the first kind (also called Stirling cycle numbers), and the 
$\Stir{n}{k}$ denotes the Stirling numbers of the second kind, respectively (see ~\cite{GraKnuPa}). Moreover, we denote with $\fallfak{x}{\ell}$ the $\ell$-th falling factorial, $\fallfak{x}\ell=x(x-1)\dots (x-\ell+1)$, $\ell\ge 0$, with $\fallfak{x}0=1$.

\section{The distributional equation}
We consider the urn model called model $M$ of Chen and Wei~\cite{ChenWei}, specified by ball replacement matrix~\eqref{NEWmodelM}.
Let $W_n$ and $B_n$ denote the random variables counting the numbers of white balls and black balls after 
$n$ draws, $n\ge 0$. We assume that the initial numbers $B_0$ of black balls and $W_0$ of white balls satisfy $B_0>0$ and
$W_0>0$, and that the total number of balls at time zero $T_0$ satisfies $T_0=B_0+W_0\ge m$,
in order to avoid any degenerate cases. Since the urn model is balanced, regardless of the inspected color combination a total of $m\cdot c$ balls 
are added to the urn at every discrete time step, the total number $T_n$ of balls after $m$ draws is a deterministic quantity, and given by 

\begin{equation}
\label{NewUrnDistrEqn1} 
T_n=T_0 + n m c = W_n+B_n,\quad n\ge 0.
\end{equation}
 
We are interested in the random variable $W_n$, counting the number of white balls contained in the urn after $n$ draws, $n\ge 0$. 
The starting point of our considerations is the distributional equation,

\begin{equation}
\label{NewUrnDistrEqn2} 
W_{n}\law W_{n-1}  + \sum_{k=0}^{m} k\cdot c \cdot \mathbb{I}_{n}(W^kB^{m-k}),
\end{equation}

which says that the number of white balls after $n$ draws can be decomposed as the number of white balls after $n-1$ draws, plus the additional balls 
added when the colors of the $n$-th draw have been inspected, $n\ge 1$. Here the random variables $\mathbb{I}_{n}(W^kB^{m-k})$, $0\le k\le m$ denote the indicators
of drawing $k$ white balls and $m-k$ black balls from the urn at the $n$-th draw, $n\ge 1$. 
Let $\mathcal{F}_{n-1}$ denote the $\sigma$-field generated by the first $n-1$ draws. By~\eqref{NewUrnDistrEqn1}, we have $B_{n-1}=T_{n-1}-W_{n-1}$, and further

\begin{equation}
\label{NewUrnDistrEqn3} 
\P\{\mathbb{I}_{n}(W^kB^{m-k})=1|\mathcal{F}_{n-1}\}=\frac{\binom{W_{n-1}}{k}\binom{B_{n-1}}{m-k} }{\binom{T_{n-1}}{m}}=\frac{\binom{W_{n-1}}{k}\binom{T_{n-1}-W_{n-1}}{m-k} }{\binom{T_{n-1}}{m}},
\end{equation}

$0\le k \le m$, $n\ge 1$. In order to study the moments $\E(W_n^s)$, $s\ge 1$, of the random variable $W_n$, we have to derive first a distributional 
equation for $W_n^s$. In order to do so, we take the $s$-th power in~\eqref{NewUrnDistrEqn2},
and use the fact that indicator variables are mutually exclusive. We obtain the distributional 
equation, valid for $n\ge 1$, $s\ge 1$:

\begin{equation}
\label{NewUrnDistrEqn4} 
W_{n}^s\law W_{n-1}^s  + \sum_{\ell=1}^{s}\binom{s}{\ell}W_{n-1}^{s-\ell}c^\ell\sum_{k=1}^{m} k^\ell  \cdot \mathbb{I}_{n}(W^kB^{m-k}).
\end{equation}

\section{Results for the moments}

\begin{theorem}
\label{NewUrnThe1}
The expected value of the random variable $W_n$, counting the numbers of white balls after 
$n$ draws, is given by $\E(W_{n})= \frac{W_0}{T_0}(nmc + T_0)$,
and the variance $\V(W_{n})=\E(W_{n}^2)-\E(W_{n})^2$ is determined via the second moment


\[
\E(W_n^2)=\frac{\binom{n-1+\lambda_{1}}{n}\binom{n-1+\lambda_{2}}{n}}{\binom{n-1+\frac{T_0}{mc}}{n}\binom{n-1+\frac{T_0-1}{mc}}{n}}
\bigg(W_0^2+\frac{W_0c^2m}{T_0}\sum_{\ell=0}^{n-1}\frac{\ell+\frac{T_0-m}{mc}}{\ell+\frac{T_0-1}{mc}}\frac{\binom{\ell+\frac{T_0}{mc}}{\ell+1}\binom{\ell+\frac{T_0-1}{mc}}{\ell+1}}{\binom{\ell+\lambda_{1}}{\ell+1}\binom{\ell+\lambda_{2}}{\ell+1}}\bigg),
\]

where the values $\lambda_1$, $\lambda_2$ are given by

\[
\lambda_{1,2}=\frac{-\frac12+mc+T_0\pm\frac12\sqrt{1+4mc(1+c)}}{mc}.
\]

\end{theorem}

Concerning higher moments, we obtain the following recursive characterization.

\begin{theorem}
\label{NewUrnThe2}
The $s$-th moment $\E(W_n^{s})$ is for $s\ge 1$ given by

\[
\E(W_{n}^s)=\bigg(\prod_{j=0}^{n-1}\alpha_{j,s}\bigg)\cdot \Big(W_0^s+  \sum_{\ell=0}^{n-1}\frac{\beta_{\ell,s}}{\prod_{j=0}^{\ell}\alpha_{j,s}}\Big),
\]

where the quantities $\alpha_{n,s},\beta_{n,s}$ are defined as
\begin{equation}
\label{NewUrnRec2}
\alpha_{n,s}= \sum_{\ell=0}^{s}c^\ell \frac{\binom{s}{\ell}\binom{m}\ell}{\binom{T_n}{\ell}},
\end{equation}

\begin{equation}
\label{NewUrnRec3}
\beta_{n,s} = \sum_{i=2}^{s}\E(W_n^{s+1-i})\sum_{\ell=i}^s\binom{s}{\ell}c^\ell 
\sum_{j=\ell+1-i}^\ell (-1)^{j+i-1-\ell}  \frac{\Stir{\ell}j\stir{j}{\ell+1-i}\binom{m}j}{\binom{T_n}{j}}.
\end{equation}

Here $\stir{n}{k}$ denotes an (unsigned) Stirling number of the first kind, and $\Stir{n}{k}$ denotes a Stirling number of the second kind.
\end{theorem}

The result above enables us to obtain in principle arbitrary high moments, 
since the moment $\E(W_n^s)$ can be expressed terms of the moments $\E(W_{h}^{r})$, with $0\le h\le n-1$, $0\le r\le s-1$. Note that vast simplifications occur in the case $m=1$, which we did not observe for $m\ge 2$. As mentioned in the introduction the case $m=1$ is very well understood; therefore we do not go into details. 
Moreover, we obtain the following results concerning the moments of the normalized random variable $W_n/n$.

\begin{coroll}
\label{NewUrnThe3}
The limits $\lim_{n\to\infty}\E(W_n^{s}/n^s)$ exist; in particular, we obtain for the expected value
$\lim_{n\to\infty}\frac{\E(W_{n})}{n}=\frac{W_0mc}{T_0}$, and for the normalized second moment
\[
\lim_{n\to\infty}\frac{\E(W_{n}^2)}{n^2}=\frac{\Gamma(\frac{T_0}{mc})\Gamma(\frac{T_0-1}{mc})}{\Gamma(\lambda_1)\Gamma(\lambda_2)}\bigg(W_0^2
+\frac{W_0c^2m}{T_0}\sum_{\ell=0}^{\infty}\frac{\ell+\frac{T_0-m}{mc}}{\ell+\frac{T_0-1}{mc}}\frac{\binom{\ell+\frac{T_0}{mc}}{\ell+1}\binom{\ell+\frac{T_0-1}{mc}}{\ell+1}}{\binom{\ell+\lambda_{1}}{\ell+1}\binom{\ell+\lambda_{2}}{\ell+1}}\bigg).
\]

Moreover, for arbitrary $s\ge 1$ the limit of the normalized $s$-th moment can be expressed in terms of an infinite sum 

\[
\lim_{n\to\infty}\frac{\E(W_{n}^s)}{n^s}=\frac{\prod_{\ell=1}^{s}\Gamma(\frac{T_0+1-\ell}{mc})}{\prod_{\ell=1}^{s}\Gamma(\lambda_{\ell,s})}
\Big(W_0^s+  \sum_{\ell=0}^{\infty}\frac{\beta_{\ell,s}}{\prod_{j=0}^{\ell}\alpha_{j,s}}\Big)
\]

with $\beta_{\ell,s}$ as defined above in Theorem~\ref{NewUrnThe2}, involving the moments of the form $\E(W_{h}^{r})$, with $0\le h<\infty$ and, $0\le r\le s-1$.
Here the $\lambda_{\ell,s}$ denote the roots (times by minus one) of the equation

\[
\frac{s!}{(mc)^s}\sum_{\ell=0}^{s}c^\ell\binom{m}\ell\binom{x\cdot mc+T_0-\ell}{s-\ell}=\prod_{\ell=1}^{s}(x+\lambda_{\ell,s}).
\]

\end{coroll}

\begin{remark}
As mentioned in the introduction Chen and Wei~\cite{ChenWei} proved almost sure convergence of $W_n$, 
namely $\mathcal{W}_n=\frac{W_n}{T_n}=\frac{W_n}{nmc+T_0}\xrightarrow{\text{a.~s.}}\mathcal{W}_{\infty}$. Hence, we get the expectation and the variance of $\mathcal{W}_{\infty}$ by the relations  $\E(\mathcal{W}_{\infty})=\lim_{n\to\infty}\frac{\E(W_{n})}{nmc}$, and $\V(\mathcal{W}_{\infty})=\lim_{n\to\infty}\frac{\V(W_{n})}{(nmc)^2}$. 
Our results show that the distribution of $\mathcal{W}_{\infty}$ is not an ordinary beta law, in contrast to the case $m=1$. Note that the moments of $W_n$ do not grow very fast, since they satisfy the trivial bounds $\E(W_n^s)\le (nmc+T_0)^s$. Hence, by Carleman's condition
the limit law $W_{\infty}$ 
is uniquely determined by its moments $\E(\mathcal{W}_{\infty}^s)$, which are given by
$\E(\mathcal{W}_{\infty}^s)=\lim_{n\to\infty}\frac{\E(W_n^s)}{(nmc)^s}$.
\end{remark}

\section{Determining the structure of the moments}

In order to prove Theorems~\ref{NewUrnThe1} and~\ref{NewUrnThe2}, we need the following result. 

\begin{lemma}
\label{NewUrnLem1}
The moments $\E(W_n^{s})$ satisfy the recurrence relation 

\begin{equation}
\label{NewUrnRec1}
\E(W_{n}^s)= \alpha_{n-1,s}\cdot\E(W_{n-1}^s)  +  \beta_{n-1,s},\quad n,s\ge 1,
\end{equation}

where the quantities $\alpha_{n,s}$ and $\beta_{n,s}$ are as defined in Theorem~\ref{NewUrnThe2}.
\end{lemma}

Our starting point for the proof of Lemma~\ref{NewUrnLem1} is the distributional equation for $W_n^s$~\eqref{NewUrnDistrEqn4}, and we take the conditional expectation 
with respect to $\mathcal{F}_{n-1}$. This leads to

\begin{equation*}
\E(W_{n}^s|\mathcal{F}_{n-1})= W_{n-1}^s  +  \sum_{\ell=1}^{s}\binom{s}{\ell}W_{n-1}^{s-\ell}c^\ell\sum_{k=1}^{m} k^\ell 
 \cdot \frac{\binom{W_{n-1}}{k}\binom{T_{n-1}-W_{n-1}}{m-k} }{\binom{T_{n-1}}{m}},
\end{equation*}

where we have used~\eqref{NewUrnDistrEqn3} and~\eqref{NewUrnDistrEqn1}.
In order to simplify the stated expression we have to use some combinatorial identities. 
We convert the ordinary powers of $k$ into falling factorials using the Stirling numbers of the second kind,

\begin{equation*}
x^m= \sum_{k=1}^{m}\Stir{m}{k}\fallfak{x}{k},
\end{equation*}

where $\fallfak{x}{k}=x(x-1)(x-2)\dots(x-(k-1))$, $m\ge 1$. We have 

\[
\sum_{k=1}^{m} k^\ell \binom{W_{n-1}}{k}=
\sum_{k=1}^{m} \sum_{j=1}^{\ell}\Stir{\ell}{j}\fallfak{k}{j}\binom{W_{n-1}}{k}
=\sum_{j=1}^{\ell} \sum_{k=j}^{m}\Stir{\ell}{j}\fallfak{W_{n-1}}{j}\binom{W_{n-1}-j}{k-j},
\]

and consequently obtain 

\begin{equation*}
\sum_{k=1}^{m} k^\ell \binom{W_{n-1}}{k}\binom{T_{n-1}-W_{n-1}}{m-k} 
=  \sum_{j=1}^{\ell} \Stir{\ell}{j}\fallfak{W_{n-1}}{j}\sum_{k=0}^{m-j}\binom{W_{n-1}-j}{k}\binom{T_{n-1}-W_{n-1}}{m-j-k}.
\end{equation*}

Next we use the Vandermonde convolution formula

\begin{equation*}
\sum_{k=0}^{n}\binom{r}{k}\binom{s}{n-k}=\binom{r+s}{n},
\end{equation*}

in order to obtain the expression

\begin{equation*}
 \sum_{k=1}^{m} k^\ell \binom{W_{n-1}}{k}\binom{T_{n-1}-W_{n-1}}{m-k} 
= \sum_{j=1}^{\ell} \Stir{\ell}{j}\fallfak{W_{n-1}}{j}\binom{T_{n-1}-j}{m-j}.
\end{equation*}

This leads to 

\[
\sum_{k=1}^{m} k^\ell 
 \cdot \frac{\binom{W_{n-1}}{k}\binom{T_{n-1}-W_{n-1}}{m-k} }{\binom{T_{n-1}}{m}}
 =\sum_{j=1}^{\ell} \Stir{\ell}{j}\fallfak{W_{n-1}}{j}\frac{\binom{T_{n-1}-j}{m-j}}{\binom{T_{n-1}}{m}}
 =\sum_{j=1}^{\ell} \Stir{\ell}{j}\fallfak{W_{n-1}}{j}\frac{\binom{m}{j}}{\binom{T_{n-1}}{j}}.
\]

Next we convert the falling factorials into ordinary powers, and obtain 

\[
\fallfak{W_{n-1}}{j}= \sum_{i=1}^{j}\stir{j}i (-1)^{j-i}W_{n-1}^i.
\]

Hence, we have 

\begin{equation*}
\begin{split}
\sum_{j=1}^{\ell} \Stir{\ell}{j}\fallfak{W_{n-1}}{j}\frac{\binom{m}{j}}{\binom{T_{n-1}}{j}}
&=\sum_{j=1}^{\ell} \Stir{\ell}{j}\frac{\binom{m}{j}}{\binom{T_{n-1}}{j}}\sum_{i=1}^{j}\stir{j}i (-1)^{j-i}W_{n-1}^i\\
&=\sum_{i=1}^{\ell}W_{n-1}^i \sum_{j=i}^{\ell}(-1)^{j-i} \frac{\stir{j}i \Stir{\ell}{j}\binom{m}{j}}{\binom{T_{n-1}}{j}}\\
&=\sum_{i=1}^{\ell}W_{n-1}^{\ell+1-i} \sum_{j=\ell+1-i}^{\ell} (-1)^{j+i-\ell-1}\frac{\stir{j}{\ell+1-i}\Stir{\ell}{j}\binom{m}{j}}{\binom{T_{n-1}}{j}}.
\end{split}
\end{equation*}
Note that the result above is an, explicit expression for the $\ell$-th moment of a hypergeometric distributed random variable with parameter $W_{n-1}$, $T_{n-1}$, and $m$.
This implies that

\begin{equation*}
\begin{split}
\E(W_{n}^s|\mathcal{F}_{n-1})&= W_{n-1}^s  +  \sum_{\ell=1}^{s}\binom{s}{\ell}W_{n-1}^{s-\ell}c^\ell\sum_{k=1}^{m} k^\ell 
 \cdot \frac{\binom{W_{n-1}}{k}\binom{T_{n-1}-W_{n-1}}{m-k} }{\binom{T_{n-1}}{m}}\\
 &=W_{n-1}^s  +  \sum_{\ell=1}^{s}\binom{s}{\ell}c^\ell\sum_{i=1}^{\ell}W_{n-1}^{s+1-i} \sum_{j=\ell+1-i}^{\ell} (-1)^{j+i-\ell-1}\frac{\stir{j}{\ell+1-i}\Stir{\ell}{j}\binom{m}{j}}{\binom{T_{n-1}}{j}} .
\end{split}
\end{equation*}

The stated result now easily follows by taking the expectation on both sides.

\begin{remark}
Note that in the case $c=1$ a simpler expression exists for the factorial moments $\E(\fallfak{W_n}{s})=\E(W_{n}(W_n-1)\dots(W_n-s+1))$ of $W_n$, and consequently 
also for the ordinary moments of $W_n$. 

\begin{equation*}
\E(\fallfak{W_n}{s})=\E(\fallfak{W_{n-1}}{s})\cdot \sum_{\ell=0}^{s} \frac{\binom{s}{\ell}\binom{m}\ell}{\binom{T_n}{\ell}}
+ \sum_{i=1}^{s}i!\E(\fallfak{W_{n-1}}{s-i})\sum_{\ell=i}^s\frac{\binom{s}\ell \binom{m}\ell \binom{s-\ell}i \binom{\ell}i}{\binom{T_{n-1}}\ell}.
\end{equation*}

\end{remark}

Next we use Lemma~\ref{NewUrnLem1} in order to prove Theorem~\ref{NewUrnThe2}; we have

\[
\E(W_{n}^s)= \alpha_{n-1,s}\cdot\E(W_{n-1}^s)  +  \beta_{n-1,s},\quad n\ge 1,s\ge 1.
\]

Let $e_{n,s}$ be defined as 

\[
e_{n,s}=\frac{\E(W_{n}^s)}{\prod_{j=0}^{n-1}\alpha_{j,s}},\quad n\ge 0, \quad\text{with}\quad
e_{0,s}=\E(W_{0}^s)=W_0^s.
\]

We have

\[
e_{n,s}= e_{n-1,s}  +  \frac{\beta_{n-1,s}}{\prod_{j=0}^{n-1}\alpha_{j,s}},\quad n\ge 1,
\]

and consequently, 

\[
e_{n,s}=\frac{\E(W_{n}^s)}{\prod_{j=0}^{n-1}\alpha_{j,s}}=e_{0,s}+\sum_{\ell=0}^{n-1}\frac{\beta_{\ell,s}}{\prod_{j=0}^{\ell}\alpha_{j,s}}
=W_0^s+\sum_{\ell=0}^{n-1}\frac{\beta_{\ell,s}}{\prod_{j=0}^{\ell}\alpha_{j,s}},
\]

which implies the stated result.

In order to obtain the result for the expected value and the variance (second moment), as stated in Theorem~\ref{NewUrnThe1}, 
we proceed as follows. In the case of the expected value we observe that $\beta_{j,1}=0$, and that 

\begin{equation*}
\alpha_{j,1}=1+\frac{cm}{T_{j}}= \frac{T_{j}+cm}{T_{n-1}}=\frac{T_{j+1}}{T_{j}},
\end{equation*}

since $T_{j+1}=T_{j}+mc$, the total number of balls contained in the urn increases by $mc$ after each draw. 
Consequently, 

\[
\prod_{j=0}^{n-1}\alpha_{j,1}=\prod_{j=0}^{n-1}\frac{T_{j+1}}{T_{j}}=\frac{T_n}{T_0}
\]

and the stated result follows.

\begin{remark}
As already mentioned before $\mathcal{W}_n=W_n/T_n\xrightarrow{\text{a.~s.}}\mathcal{W}_{\infty}$. 
This can easily be seen as follows: we readily note that for $s=1$ it holds

\[
\E(W_{n}|\mathcal{F}_{n-1})=W_{n-1}\frac{T_{n}}{T_{n-1}},\quad n\ge 1.
\]

Thus, $\mathcal{W}_n=W_n/T_n$ is a positive martingale with respect to the filtration $\mathcal{F}_n$, as previously observed
in~\cite{ChenWei}, which directly leads to the proof of the almost sure convergence of $\mathcal{W}_n$.
\end{remark}

In order to obtain the variance $\V(W_n)=\E(W_n^2)-\E(W_n)^2$, we study the second moment $\E(W_n^2)$, given by

\[
\E(W_{n}^2)=\bigg(\prod_{j=0}^{n-1}\alpha_{j,2}\bigg)\cdot \Big(W_0^2+  \sum_{\ell=0}^{n-1}\frac{\beta_{\ell,2}}{\prod_{j=0}^{\ell}\alpha_{j,2}}\Big).
\]

The value $\alpha_{n,2}$ is given by

\begin{equation*}
\begin{split}
\alpha_{n,2}&
=1+\frac{2cm}{T_n}+\frac{c^2\binom{m}{2}}{\binom{T_n}{2}}=
\frac{T_{n+2}(T_{n}-1)+c^2m(m-1)}{T_n(T_{n}-1)}\\
&=\frac{(n+2+\frac{T_0}{cm})(n+\frac{T_0-1}{cm})+1-\frac1m}{(n+\frac{T_0}{cm})(n+\frac{T_0-1}{cm})}.
\end{split}
\end{equation*}

We factor the numerator of $\alpha_{n,2}$ by determining the zeros of the quadratic equation with respect to $n$ and 
get

\[
\alpha_{n,2}=\frac{(n+\lambda_1)(n+\lambda_2)}{(n+\frac{T_0}{cm})(n+\frac{T_0-1}{cm})},
\]

with $\lambda_1,\lambda_2$ as stated in Theorem~\ref{NewUrnThe1}. Concerning $\beta_{n,2} $ we have

\begin{equation*}
\beta_{n,2} = 
\E(W_n)c^2\Big(\frac{m}{T_n}-\frac{\binom{m}{2}}{\binom{T_n}{2}}\Big)
=\frac{W_0c^2m}{T_0}\Big(1-\frac{m-1}{T_n-1}\Big)=\frac{W_0c^2m}{T_0}\cdot\frac{T_n-m}{T_n-1}.
\end{equation*}

This readily leads to the stated exact result for the second moment.

\subsection{Asymptotic expansions}
We finally turn to the proof of Corollary~\ref{NewUrnThe3}. 
In order to get the results for the limit $\lim_{n\to\infty}\E(W_n^2)/n^2$
we proceed as follows: first we write $\prod_{j=0}^{n-1}\alpha_{j,2}$ in terms of Gamma-functions,

\[
\prod_{j=0}^{n-1}\alpha_{j,2}=
\frac{\binom{n-1+\lambda_{1}}{n}\binom{n-1+\lambda_{2}}{n}}{\binom{n-1+\frac{T_0}{mc}}{n}\binom{n-1+\frac{T_0-1}{mc}}{n}}
=\frac{\Gamma(n+\lambda_1)\Gamma(n+\lambda_2)\Gamma(\frac{T_0}{mc})\Gamma(\frac{T_0-1}{mc})}{\Gamma(\lambda_1)\Gamma(\lambda_2)\Gamma(n+\frac{T_0}{mc})\Gamma(n+\frac{T_0-1}{mc})}.
\]

Using Stirling's formula for the Gamma function
\begin{equation}
\label{NEWStirling}
    \Gamma(z) = \Bigl(\frac{z}{e}\Bigr)^{z}\frac{\sqrt{2\pi }}{\sqrt{z}}%
    \Bigl(1+\frac{1}{12z}+\frac{1}{288z^{2}}+\mathcal{O}(\frac{1}{z^{3}})\Bigr),
\end{equation}

and the fact that $\lambda_1+\lambda_2=2+\frac{2W_0-1}{mc}$,
we obtain the asymptotic expansion

\[
\prod_{j=0}^{n-1}\alpha_{j,2}=
n^2\frac{\Gamma(\frac{T_0}{mc})\Gamma(\frac{T_0-1}{mc})}{\Gamma(\lambda_1)\Gamma(\lambda_2)}\Big(1+\mathcal{O}\big(\frac1n\big)\Big).
\]

The asymptotic expansion of $\prod_{j=0}^{n-1}\alpha_{j,2}$ also implies that the sum $\sum_{\ell=0}^{n-1}\frac{\beta_{\ell,2}}{\prod_{j=0}^{\ell}\alpha_{j,2}}$ is convergent for $n\to\infty$ 
by the comparison test, 

\[
\sum_{\ell=0}^{n-1}\frac{\beta_{\ell,2}}{\prod_{j=0}^{\ell}\alpha_{j,2}}\le K\sum_{\ell=1}^{n-1}\frac{1}{\ell^2},
\]
for some suitable constant $K>0$.

In order to provide the asymptotics of $\E(W_{n}^s)$ we proceed by induction. Before we actually do so, 
we derive an asymptotic expansion of $\prod_{j=0}^{n-1}\alpha_{j,s}$. We have

\begin{equation*}
\begin{split}
\alpha_{j,s}&=\sum_{\ell=0}^{s}c^\ell \frac{\binom{s}{\ell}\binom{m}\ell}{\binom{T_j}{\ell}}
=\frac{\sum_{\ell=0}^{s}c^\ell\binom{m}\ell \binom{T_j-\ell}{s-\ell}}{\binom{T_j}{s}}
=s!\frac{\sum_{\ell=0}^{s}c^\ell\binom{m}\ell \binom{T_j-\ell}{s-\ell}}{(mc)^s\prod_{\ell=0}^{s-1}(j+\frac{T_0-\ell}{mc})}\\
&=s!\frac{\sum_{\ell=0}^{s}c^\ell\binom{m}\ell \binom{jmc+T_0-\ell}{s-\ell}}{(mc)^s\prod_{\ell=1}^{s}(j+\frac{T_0+1-\ell}{mc})}
\end{split}
\end{equation*}

Let $\lambda_{\ell,s}$ denote the roots (times minus one) of the equation

\[
\frac{s!}{(mc)^s}\sum_{\ell=0}^{s}c^\ell\binom{m}\ell\binom{x\cdot mc+T_0-\ell}{s-\ell}=\prod_{\ell=1}^{s}(x+\lambda_{\ell,s}).
\]

We get

\[
\alpha_{j,s}
=s!\frac{\sum_{\ell=0}^{s}c^\ell\binom{m}\ell \binom{jmc+T_0-\ell}{s-\ell}}{(mc)^s\prod_{\ell=1}^{s}(j+\frac{T_0+1-\ell}{mc})}
=\frac{\prod_{\ell=1}^{s}(j+\lambda_{\ell,s})}{\prod_{\ell=1}^{s}(j+\frac{T_0+1-\ell}{mc})},
\]

and consequently 

\[
\prod_{j=0}^{n-1}\alpha_{j,s}=\prod_{\ell=1}^{s}\frac{\Gamma(n+\lambda_{\ell,s})\Gamma(\frac{T_0+1-\ell}{mc})}{\Gamma(\lambda_{\ell,s})\Gamma(n+\frac{T_0+1-\ell}{mc})}.
\]

By Stirling's formula we obtain the asymptotic expansion

\[
\prod_{j=0}^{n-1}\alpha_{j,s}=\frac{n^{\sum_{\ell=1}^{s}\lambda_{\ell,s}}}{n^{\frac{sT_0-\binom{s}{2}}{mc}}}\prod_{\ell=1}^{s}\frac{\Gamma(\frac{T_0+1-\ell}{mc})}{\Gamma(\lambda_{\ell,s})}\Big(1+\mathcal{O}\big(\frac1n\big)\Big).
\]

Let $[x^k]$ denote the extraction of coefficients operator. Since 

\begin{equation*}
\begin{split}
\sum_{\ell=1}^{s}\lambda_{\ell,s}&=[x^{s-1}]\prod_{\ell=1}^{s}(x+\lambda_{\ell,s})
=[x^{s-1}]\frac{s!}{(mc)^s}\sum_{\ell=0}^{s}c^\ell\binom{m}\ell\binom{x\cdot mc+T_0-\ell}{s-\ell}\\
&=\frac{s!}{(mc)^s}\Big(\frac{(cm)^{s-1}(sT_0-\binom{s}{2})}{s!}+\frac{cm\cdot(cm)^{s-1}}{(s-1)!}\Big)
=\frac{sT_0-\binom{s}2}{mc} + s,
\end{split}
\end{equation*}

we obtain 

\[
\prod_{j=0}^{n-1}\alpha_{j,s}=n^{s}\prod_{\ell=1}^{s}\frac{\Gamma(\frac{T_0+1-\ell}{mc})}{\Gamma(\lambda_{\ell,s})}\Big(1+\mathcal{O}\big(\frac1n\big)\Big).
\]

Now we assume that $\E(W_{n}^s)=\kappa_s\cdot n^s(1+\mathcal{O}\big(\frac1n\big)$, being true for $s=1,2$. By Theorem~\ref{NewUrnThe2} we have
 
\begin{equation*}
\begin{split}
\E(W_{n}^s)&=\bigg(\prod_{j=0}^{n-1}\alpha_{j,s}\bigg)\cdot \Big(W_0^s+  \sum_{\ell=0}^{n-1}\frac{\beta_{\ell,s}}{\prod_{j=0}^{\ell}\alpha_{j,s}}\Big)\\
&=n^{s}\prod_{\ell=1}^{s}\frac{\Gamma(\frac{T_0+1-\ell}{mc})}{\Gamma(\lambda_{\ell,s})}\Big(1+\mathcal{O}\big(\frac1n\big)\Big)
\cdot  \Big(W_0^s+  \sum_{\ell=0}^{n-1}\frac{\beta_{\ell,s}}{\prod_{j=0}^{\ell}\alpha_{j,s}}\Big).
\end{split}
\end{equation*}

By our induction hypothesis, $\E(W_{\ell}^{k})=\kappa_k\cdot \ell^k(1+\mathcal{O}\big(\frac1\ell\big)$ for sufficiently large $\ell$, $1\le k\le s-1$,
and consequently $\beta_{\ell,s}$ satisfies $\beta_{\ell,s}\le K\cdot \ell^{s-2}$, where $K$ denotes a sufficiently large constant, only depending on $s$. 
Using our result for $\prod_{j=0}^{\ell}\alpha_{j,s}$, the sum $\sum_{\ell=0}^{n-1}\frac{\beta_{\ell,s}}{\prod_{j=0}^{\ell}\alpha_{j,s}}$ 
is convergent according to the comparison test with $\sum_{\ell=1}^{\infty}\frac{1}{\ell^2}$. Hence, $\E(W_{n}^s)=\kappa_s\cdot n^s(1+\mathcal{O}\big(\frac1n\big))$, 
and we have proven our stated results.

\section{The case of three or more colors}

The urn model $M$ can be readily generalized to $r\ge 2$ different colors. As in the case of two colors $m$ balls are drawn at random from the urn, say $k_i$ balls
of color $i$, $1\le i \le r$, their colors are noted, and they are returned to the urn together with $c\cdot k_i$ balls of color $i$.
Let $\mathbf{X}_n=(X_{n,1},\dots,X_{n,r})$ 
denote the random vector counting the number of balls of type $i$ contained in the urn after $n$ draws, with 
initial values $\mathbf{X}_0=(X_{0,1},\dots,X_{0,r})$. 
The urn is again balanced, so the number of balls after $n$ draws is given by $T_n=T_0 + n m c$, with $T_0=\sum_{i=1}^{r}X_{0,i}$.
One gets the distributional equation

\begin{equation*}
\mathbf{X}_n\law \mathbf{X}_{n-1} + \sum_{\substack{k_1+\dots+ k_r=m\\k_\ell\ge 0}} c\cdot \mathbf{k} \cdot \mathbb{I}_{n}(\prod_{\ell=1}^{r}X_{\ell}^{k_\ell}),
\end{equation*}

where $\mathbf{k}=(k_1,\dots,k_r)$, and $\mathbb{I}_{n}(\prod_{\ell=1}^{r}X_{\ell}^{k_\ell})$ denote the indicators of drawing $k_\ell$ balls of colour $\ell$, $1\le\ell\le r$ at the $n$-th draw. Note that each individual random variable $X_{n,\ell}$ satisfies $X_{n,\ell}\law W_{n}$, where $W_{n}$ denotes the previously considered random variable from the two color case. The distributional equation above can be used to study the mixed moments of $X_{n,\ell}$ similar to the results of Theorems~\ref{NewUrnThe1},~\ref{NewUrnThe2}, and Corollary~\ref{NewUrnThe3}.
We refrain from going into details since the resulting expressions get very involved; we only mention our findings for the 
covariance of two different colors $i,j$, with $r\ge 3$ and $1\le i<j\le m$.

\[
\text{Cov}(X_{n,i},X_{n,j})=\frac{\binom{n-1+\lambda_{1}}{n}\binom{n-1+\lambda_{2}}{n}}{\binom{n-1+\frac{T_0}{mc}}{n}\binom{n-1+\frac{T_0-1}{mc}}{n}}X_{0,i}X_{0,j}
- \frac{(mc)^2(n+\frac{T_0}{mc})^2}{T_0^2}X_{0,i}X_{0,j}.
\]

Moreover, in the limit we obtain 

\[
\lim_{n\to\infty}\text{Cov}(\frac{X_{n,i}}{n},\frac{X_{n,j}}{n})=X_{0,i}X_{0,j}\Big(\frac{\Gamma(\frac{T_0}{mc})\Gamma(\frac{T_0-1}{mc})}{\Gamma(\lambda_1)\Gamma(\lambda_2)}
-\frac{(mc)^2}{T_0^2}\Big),
\]

with $\lambda_1,\lambda_2$ as given in Theorem~\ref{NewUrnThe1}.

\section{Urn model R- drawing with replacement}

We consider another urn model which we call model $R$. It can be considered as a variant of model $M$. 
The dynamics of model $R$ are identical to model $M$: $m\ge 1$ balls are drawn from the urn containing balls of colors black and white. Their colors are inspected, and they are returned to the
urn. According to the observed colors we add new balls: if $\ell$ black balls and $m-\ell$ white balls have been observed, 
we add  $c(m - \ell)$ white balls and $\ell$ black balls; the ball replacement matrix coincides with~\eqref{NEWmodelM}. 
The main difference in model $R$ is that the sampling of the $m$ balls occurs \emph{with replacement}, i.e.~ the $m$ balls
are drawn one by one from the urn, the colors are observed, and then put back into the urn. Hence, the distribution of the number of white balls
in the sample of size $m$ follows a binomial distribution instead of a hypergeometric distribution in model $M$. 
The distributional equation for the number of white balls $W_n$ after $n$ draws is identical to~\eqref{NewUrnDistrEqn2},

\[
W_{n}\law W_{n-1}  + \sum_{k=0}^{m} k\cdot c \cdot \mathbb{I}_{n}(W^kB^{m-k}),
\]

but the distribution of the indicator variables $\mathbb{I}_{n}(W^kB^{m-k})$ changes to 

\begin{equation}
\label{NewUrnModelReqn1} 
\P\{\mathbb{I}_{n}(W^kB^{m-k})=1|\mathcal{F}_{n-1}\}=\frac{\binom{m}kW_{n-1}^kB_{n-1}^{m-k} }{T_{n-1}^m}
=\frac{\binom{m}kW_{n-1}^k(T_{n-1}-W_{n-1})^{m-k} }{T_{n-1}^m},
\end{equation}
$0\le k \le m$, $n\ge 1$. 

Concerning a limit law for the number of white balls $W_n$ after $n$ draws, we can extend the results for model $M$ of Chen and Wei~\cite{ChenWei} to model $R$.

\begin{theorem}
\label{NewUrnModelRthe1}
The random variable $\mathcal{W}_n=W_n/T_n$ is a positive martingale with respect to the natural filtration $\mathcal{F}_n$, 
$\E(\mathcal{W}_n|\mathcal{F}_{n-1})=\mathcal{W}_{n-1}$. Consequently,

\[
\mathcal{W}_n\xrightarrow{a.s.}\mathcal{W}_{\infty}.
\]

Furthermore, for fixed $W_0$, $B_0$, $c$ and $m$, $\mathcal{W}_{\infty}$ is absolutely continuous. 
\end{theorem}

Moreover, the moments of $W_n$ satisfy recurrence relations similar to model $M$.

\begin{theorem}
\label{NewUrnModelRthe2}
The expected value is given by $\E(W_n)=\frac{W_0}{T_0}(nmc+T_0)$ and the variance $\V(W_n)=\E(W_n^2)-\E(W_n)^2$ in terms of second moment,

\[
\E(W_n^2)=\frac{\binom{n-1+\mu_{1}}{n}\binom{n-1+\mu_{2}}{n}}{\binom{n-1+\frac{T_0}{mc}}{n}^2}
\bigg(W_0^2+\frac{W_0c^2m}{T_0}\sum_{\ell=0}^{n-1}\frac{\binom{\ell+\frac{T_0}{mc}}{\ell+1}^2}{\binom{\ell+\mu_{1}}{\ell+1}\binom{\ell+\mu_{2}}{\ell+1}}\bigg)
\]

where the values $\mu_1$, $\mu_2$ are given by $\mu_{1,2}=\frac{T_0+mc\pm c\sqrt{m}}{mc}$.
The $s$-th moment satisfies the recurrence relation

\[
\E(W_{n}^s)=\bigg(\prod_{j=0}^{n-1}\gamma_{j,s}\bigg)\cdot \Big(W_0^s+  \sum_{\ell=0}^{n-1}\frac{\delta_{\ell,s}}{\prod_{j=0}^{\ell}\gamma_{j,s}}\Big),
\]

where the quantities $\gamma_{n,s},\delta_{n,s}$ are defined as
\begin{equation*}
\gamma_{n,s}= \sum_{\ell=0}^{s}c^\ell \frac{\binom{s}{\ell}\fallfak{m}{\ell}}{T_n^\ell},
\quad \delta_{n,s} = \sum_{j=2}^{s}\E(W_n^{s+1-j})\sum_{\ell=j}^s\binom{s}{\ell}C^\ell 
\frac{\Stir{\ell}{\ell+1-j}\fallfak{m}{\ell+1-j}}{T_n^{\ell+1-j}}.
\end{equation*}

\end{theorem}

In the following we present recurrence relations for the limits of the moments of the normalized random variable $W_n/n$, 
with $\lim_{n\to\infty}\E(W_n^s/n^s)=(mc)^s\E(\mathcal{W}_\infty^s)$.

\begin{coroll}
\label{NewUrnModelRthe3}
The limits of the normalized moments $\E(W_n^s/n^s)$ exist, , 
and satisfy

\[
\lim_{n\to\infty}\frac{\E(W_{n}^s)}{n^s}=\frac{\Gamma(\frac{T_0}{mc})^s}{\prod_{\ell=1}^{s}\Gamma(\mu_{\ell,s})}
\Big(W_0^s+  \sum_{\ell=0}^{\infty}\frac{\delta_{\ell,s}}{\prod_{j=0}^{\ell}\gamma_{j,s}}\Big)
\]

with $\gamma_{j,s},\delta_{\ell,s}$ as defined above. Here the $\mu_{\ell,s}$ denote the roots (times minus one) of the equation

\[
\frac{1}{(mc)^s}\sum_{\ell=0}^{s}\binom{s}\ell\binom{m}\ell c^\ell\ell!(x\cdot mc+T_0)^{s-\ell}=\prod_{\ell=1}^{s}(x+\mu_{\ell,s}).
\]

\end{coroll}

In the following we first sketch the proofs of Theorem~\ref{NewUrnModelRthe2} and Corollary~\ref{NewUrnModelRthe3}.
Since the proofs are very similar to the proofs for model $M$, we will be very brief.
Then we discuss the proof of Theorem~\ref{NewUrnModelRthe1}.

\subsection{The structure of the moments}
Our starting point is again the distributional equation for $W_n$, which leads to a distributional equation for $W_n^s$. 
We take the conditional expectation with respect to $\mathcal{F}_{n-1}$, and obtain

\begin{equation*}
\E(W_{n}^s|\mathcal{F}_{n-1})= W_{n-1}^s  +  \sum_{\ell=1}^{s}\binom{s}{\ell}W_{n-1}^{s-\ell}c^\ell\sum_{k=1}^{m} k^\ell \binom{m}k
 \frac{W_{n-1}^k(T_{n-1}-W_{n-1})^{m-k} }{T_{n-1}^m}.
\end{equation*}

The sum appearing on the right hand side is the $\ell$-th moment of a binomial distributed random variable with parameters $m$ and $W_{n-1}/T_{n-1}$. 
We get the result 

\[
\sum_{k=1}^{m} k^\ell \binom{m}k
\frac{W_{n-1}^k(T_{n-1}-W_{n-1})^{m-k} }{T_{n-1}^m}=
 \sum_{j=1}^{\ell}\Stir{\ell}j \fallfak{m}{j}\frac{W_{n-1}^j}{T_{n-1}^j}.
\]

This implies that 

\[
\E(W_{n}^s|\mathcal{F}_{n-1})= \gamma_{n-1,s}W_{n-1}+\delta_{n-1,s},
\]
with $\gamma_{n,s}$, $\delta_{n,s}$ as stated in Theorem~\ref{NewUrnModelRthe2}.
This recurrence relation can be solved similar to the proof of Theorem~\ref{NewUrnThe2} for model $M$, and the stated result for $\E(W_n^s)$ follows.
Moreover, we easily obtain the stated formula for the expected value, and also the result for the second moment of $W_n$.
The results for the higher moments can be obtained similar to model $M$. Concerning the asymptotic expansions one can proceed similar to the proof of Corollary~\ref{NewUrnThe3}; we omit the details.
\smallskip

\subsection{Martingales and absolute continuouity}
Since $\delta_{n-1,1}=0$ we obtain for $\mathcal{W}_{n}=W_{n}/T_n$

\[
\E(\mathcal{W}_{n}|\mathcal{F}_{n-1})=\mathcal{W}_{n-1}.
\]

Hence, $\mathcal{W}_n$ is a martingale. Since it is a positive martingale, it converges almost surely to a limit $\mathcal{W}_{\infty}$. 
Note that, one has $\mathcal{W}_n \sim\frac{W_n}{cmn}$ provided $W_n\to\infty$, and thus we can obtain the moments of $\mathcal{W}_{\infty}$ 
via the moments of $\frac{W_n}{n}$.

\smallskip

Concerning the absolute continuouity of the distribution of $\mathcal{W}_\infty$ one can adapt the argumentation of Chen and Wei~\cite{ChenWei}. 
For the convenience of the reader we outline the main steps, quote the main results of~\cite{ChenWei}, and present the new ingredient of the proof for model $R$ in Lemma~\ref{NewUrnLemmaEinfach}. Let $(\Omega, \mathcal{F} , \P)$ be the probability space. In order to prove that the distribution of $\mathcal{W}_{\infty}$ has a density, one introduces a sequence of events $(\Omega_{\ell})_{\ell\ge 1}$ such that $\Omega_{\ell}\subset \Omega_{\ell+1}$, $\P\{\cup_{\ell=1}^{\infty}\Omega_{\ell}\}=1$. Then, $\mathcal{W}_\infty$ is restricted to $\Omega_{\ell}$ and it is shown that it has a density $f_{\ell}$. The proof is then finished by proving that $f =\lim_{\ell\to\infty}f_{\ell}$ exists and that $f$ is the density of $\mathcal{W}_\infty$. 

\begin{prop}[Chen and Wei~\cite{ChenWei}]
\label{NewUrnChenWei1}
Let$(\Omega_{\ell})_{\ell\ge 1}$be a sequence of increasing events such that $\P\{\cup_{\ell=1}^{\infty}\Omega_{\ell}\}=1$. If there
exist nonnegative Borel measurable functions $(f_\ell)_{\ell\ge 1}$ such that 
$\P(\Omega_{\ell}\cap\mathcal{W}_\infty^{-1}(B)) = \int_B f_{\ell}(x)dx$
for all Borel sets $B$, then$f =\lim_{\ell\to\infty}f_{\ell}$ exists almost everywhere, and $f$ is the density
of $\mathcal{W}_\infty$.
\end{prop}

In order to construct the sequence of events $(\Omega_{\ell})_{\ell\ge 1}$ 
one can follow~\cite{ChenWei}:

\begin{prop}
\label{NewUrnChenWei2}
For fixed $W_0$, $B_0$, $c$ and $m$ let

\[
\Omega_{\ell}=\{\omega: W_\ell(\omega)\ge cm \quad\text{and}\quad B_\ell(\omega)\ge cm\}, \quad \ell\ge 1.
\]

Then, $\Omega_{\ell}\subset \Omega_{\ell+1}$, and $\P\{\cup_{\ell=1}^{\infty}\Omega_{\ell}\}=1$.
\end{prop}

In order to show that $\mathcal{W}_\infty$ has a density by restricting $\mathcal{W}_\infty$ to $\Omega_{\ell}$, 
it suffices to show that the restriction of $\mathcal{W}_\infty$ to $\Omega_{\ell,j}=\{\omega:W_{\ell}(\omega)=j\}$ has a density for each $j$ ,
with $cm \le j \le T_{\ell-1}$. For this, the following result is needed, which is main new ingredient in our argumentation (compared to model $M$).

\begin{lemma}
\label{NewUrnLemmaEinfach}
The sum $\sum_{i=0}^{m}\P\{W_{n+1}=j+ck|W_n=j+c(k-i)\}$ satisfies

\begin{equation*}
\begin{split}
&\sum_{i=0}^{m}\P\{W_{n+1}=j+ck|W_n=j+c(k-i)\}\\
&\qquad=\frac{1}{T_n^m}\sum_{\ell=0}^{m}T_n^{\ell}\sum_{i=0}^{m-\ell}\binom{m}{i}\binom{m-i}{\ell}(j+c(k-i))^i(-j-c(k-i))^{m-i-\ell}.
\end{split}
\end{equation*}

Consequently, for a fixed $\ell$ it holds for all $n\ge \ell$, $cm\le j\le T_{\ell-1}$, and $k<m(n+1)$, and a suitably choosen constant $\kappa>0$.

\[
\sum_{i=0}^{m}\P\{W_{n+1}=j+ck|W_n=j+c(k-i)\}\le 1-\frac1n+\frac{\kappa}{n^2}.
\]

\end{lemma}

\begin{remark}
Note that the corresponding result of Chen and Wei for model $M$ (Lemma 4.1 and Lemma 4.2 in~\cite{ChenWei}) 
can be extended and largely simplified noting that

\begin{equation*}
\begin{split}
&\sum_{i=0}^{m}\binom{j+c(k-i)}{i}\binom{T_n-j-c(k-i)}{m-i}=\\
&\qquad\sum_{f=0}^{m}T_n^{f}\sum_{i=0}^{m-f}\binom{j+c(k-i)}{i}\sum_{\ell=f}^{m-i}\frac{\stir{\ell}f (-1)^{\ell-f}\binom{-j-c(k-i)}{m-i-\ell}}{\ell!}.
\end{split}
\end{equation*}
\end{remark}

\begin{proof}[Proof of Lemma~\ref{NewUrnLemmaEinfach}]
One has

\begin{equation*}
\begin{split}
&\sum_{i=0}^{m}\P\{W_{n+1}=j+ck|W_n=j+c(k-i)\}\\
&\quad=\frac{1}{T_n^{m}}\sum_{i=0}^{m}\binom{m}{i}(j+c(k-i))^i(T_n-j-c(k-i))^{m-i}\\
&\quad=\frac{1}{T_n^{m}}\sum_{i=0}^{m}\binom{m}{i}(j+c(k-i))^i\sum_{\ell=0}^{m-i}\binom{m-i}{\ell}T_n^{\ell}(-j-c(k-i))^{m-i-\ell}\\
&\quad =\frac{1}{T_n^m}\sum_{\ell=0}^{m}T_n^{\ell}\sum_{i=0}^{m-\ell}\binom{m}{i}\binom{m-i}{\ell}(j+c(k-i))^i(-j-c(k-i))^{m-i-\ell},
\end{split}
\end{equation*}

such that 

\begin{equation*}
\begin{split}
&\sum_{i=0}^{m}\P\{W_{n+1}=j+ck|W_n=j+c(k-i)\}
=1-\frac{(m-1)(j+ck)+c}{T_n}+\mathcal{O}(\frac1{n^2}).
\end{split}
\end{equation*}

Using $T_n=nmc+T_0$ we get further 

\[
\frac{(m-1)(j+ck)+c}{T_n}=\frac{(m-1)(j+ck)+c}{nmc}+\mathcal{O}(\frac1{n^2}).
\]

Since $\frac{(m-1)(j+ck)+c}{mc}\ge 1$, for $0\le k\le m(n+1)$, $W_0+cm\le j\le T_{\ell-1}$, the stated result follows.
\end{proof}

The proof of Theorem~\ref{NewUrnModelRthe1} can be finished by combining the results of Propositions~\ref{NewUrnChenWei1},~\ref{NewUrnChenWei2}
and Lemma~\ref{NewUrnLemmaEinfach}; it is identical to the proof in~\cite{ChenWei} (Proof of Theorem 4.2) and therefore the details are omitted.

\section{Outlook}
An interesting line of research would be to extend the results for models $M$ and $R$ to non-balanced urns with 
ball replacement matrix given by

\begin{equation}
\left[
	\begin{array}{cc}
	ma&	0	\\
	(m-1)a	& b\\
	\hdots&\hdots\\
	a	& (m-1)b\\
	0&	mb	\\
	\end{array}
	\right],
\end{equation}

with $a,b\in\N$. In contrast to the balanced versions it seems much more difficult to obtain closed form expressions for the moments of $W_n$.

\smallskip

One can analyze a generalized Friedman's urn, 
with ball replacement matrix 

\begin{equation}
\label{NEWmodelFriedman}
	M=\left[
	\begin{array}{cc}
	0&	mc	\\
	c	& (m-1)c\\
	\hdots&\hdots\\
	(m-1)c	& c\\
	mc&	0	\\
	\end{array}
	\right],
\end{equation}

and parameters $c\in\N$, and $m\ge 1$. If, say, $m - \ell$ white balls and $\ell$ black balls have been drawn from the urn, $0\le\ell\le m$, then the drawn balls are returned to the urn together with addition $c\ell$ white balls and $c(m-\ell)$ black balls. 
It is possible to look at drawing the $m$ balls without replacement, which we call model $FM$, 
or to look at drawing the $m$ balls with replacement, or simply model $FR$.
In any case, the distributional equation for the number of white balls $W_n$ after $n$ draws is given by

\[
W_{n}\law W_{n-1}  + \sum_{k=0}^{m} (m-k)\cdot c \cdot \mathbb{I}_{n}(W^kB^{m-k}),
\]

with

\[
\P\{\mathbb{I}_{n}(W^kB^{m-k})=1|\mathcal{F}_{n-1}\}=
\begin{cases}
\displaystyle{\frac{\binom{W_{n-1}}{k}\binom{T_{n-1}-W_{n-1}}{m-k} }{\binom{T_{n-1}}{m}}},&\text{model }FM\\[0.4cm]
\displaystyle{\frac{\binom{m}k W_{n-1}^k(T_{n-1}-W_{n-1})^{m-k} }{T_{n-1}^m}}, &\text{model }FR\\
\end{cases}
\]

$0\le k \le m$, $n\ge 1$. We can obtain the moments of $W_n$ for both models $FM$ and $FR$, similar 
to our previous results for models $M$ and $R$. In particular, the expectation and the variance can be obtained; we get for example the result

\[
 \E(W_n) =
 \begin{cases}
  \displaystyle{\frac{(mc)^2\binom{n}2 + mcT_0 n + (T_0-mc)W_0}{mc(n-1)+T_0}}, & mc\neq T_0\\[0.4cm] 
  \displaystyle{\frac{mc(n+1)}{2}}, & mc=T_0,\\[0.1cm] 
  \end{cases}
\]

$n\ge 0$, being valid for both models $FM$ and $FR$. In contrast to the generalized P\'olya urn models $M$ and $R$
the simple martingale structure of $W_n/T_n$ is not present anymore. However, one can find values $\varphi_n$ and $\psi_n$
such that $\mathcal{M}_n=\varphi_n W_n+\psi_n$ is a martingale,

\[
\E(\mathcal{M}_n|\mathcal{F}_{n-1})=\mathcal{M}_{n-1},
\]

with $\varphi_n$ and $\psi_n$ being given by

\[
\varphi_n= 
\frac{T_{n-1}}{T_0},
\qquad
\psi_n=
-cm\sum_{k=0}^{n-1}\frac{T_k}{T_0}.
\]


It seems to be possible to derive a limit law using the martingale central limit theorem, or to apply the methods of moments in order to obtain moment convergence. 
The study of the generalized Friedman urn models $FM$ and $FR$ will be the subject of a forthcoming companion work.

\section*{Acknowledgements}
This research was started during the second author's stay at the Academia Sinica; he wants to thank the institute for its hospitality and for providing excellent working conditions, and in particular Hsien-Kuei Hwang for interesting discussions about urn models.

\end{document}